\documentclass[a4paper,12pt]{article}
\usepackage{amssymb,enumerate}
\usepackage{amsmath,amssymb}
\newtheorem{theorem}{Theorem}
\newtheorem{df}{Definition}

\newtheorem{corollary}{Corollary}

\def\bp{{\noindent\bf Proof. \ }}

\title{Two classes of Meir-Keeler contractions}

\author{{L. G\u avru\c ta, P. G\u avru\c ta, F. Khojasteh} }
\date{}

\begin{document}

\maketitle

\begin{minipage}{120mm}
\small{\bf Abstract.} {In the present paper, we prove that $\mathcal{Z}-$contractions and weakly type contractions are actually Meir-Keeler contractions.}\\



\end{minipage}

Let $(X,d)$ be a metric space and $T:X\rightarrow X$ be a mapping. $T$ is called a Banach contraction on $X$ if there exist $\lambda\in[0,1)$ so that $$d(Tx, Ty)\leq\lambda d(x,y),\quad\textrm{for all}~x,y\in X.$$

 S. Banach proved in \cite{Banach} that $T$ has a unique fixed point $u$ in $X$ and for every $x_0\in X$ the Picard sequence $\{x_n\}$, where $x_n=Tx_{n-1},$ for all $n\in\mathbb{N}$, converges to the fixed point of $T.$
 After this result, a large number of generalizations were obtained. See, for example, the book of I.A. Rus et al. \cite{Rus} and the articles \cite{Cadariu} and \cite{Rhoades}.
 \begin{df} We say that  $T$ is a Meir-Keeler contraction if given an $\varepsilon>0,$ there exists $\delta>0$ such that $$\varepsilon\leq d(x,y)<\varepsilon+\delta ~\textrm{implies}~d(Tx, Ty)<\varepsilon,$$
 for all $x,y\in X.$
 \end{df}
 In \cite{Meir} it is proved the following Theorem.
\begin{theorem}\label{thmufp}Let $(X,d)$ be a complete metric space and $T$ be a Meir-Keeler contraction. Then $T$ has a unique fixed point $u$ and the Picard sequence $x_n=Tx_{n-1}$ converges to the fixed point of $T$.
\end{theorem}

In \cite{Khojasteh}, the authors introduced a new class of contractions, called $\mathcal{Z}-$contractions.
\begin{df}\label{defsimfct}
A mapping $\zeta:[0,\infty)\times[0,\infty)\rightarrow\mathbb{R}$ is called a simulation function if it satisfies the following conditions
\begin{enumerate}[$i)$]
\item $\zeta(0,0)=0$;
\item $\zeta(t,s)<s-t$, for all $t,s>0$;
\item if $\{t_n\}$, $\{s_n\}$ are sequences in $(0,\infty)$ such that $\displaystyle \lim_{n\rightarrow\infty}t_n=\lim_{n\rightarrow\infty}s_n>0$ then $\displaystyle \limsup_{n\rightarrow\infty}\zeta(t_n,s_n)<0$
\end{enumerate}
We denote by $\mathcal{Z}$ the set of all simulation functions.
\end{df}

\begin{df}
$T$ is called a $\mathcal{Z}$-contraction if $$\zeta(d(Tx, Ty), d(x,y))\geq 0,\quad \textrm{for all}~ x,y\in X.$$
\end{df}
The authors of the paper \cite{Khojasteh} proved the following Theorem.
\begin{theorem}Let $(X,d)$ be a complete metric space and $T$ a $\mathcal{Z}$-contraction. Then $T$ has a unique fixed point $u$ and the Picard sequence $x_n=Tx_{n-1}$ converges to $u.$
\end{theorem}
We prove that the above result is follows from the following Theorem.
\begin{theorem} A $\mathcal{Z}$-contraction is a Meir-Keeler contraction contraction.
\end{theorem}
\bp Let be $T:X\rightarrow X$ a $\mathcal{Z}$-contraction. We suppose that $T$ is not a Meir-Keeler. Then there is $\varepsilon_0>0$ so that for any $\delta>0$ there is $x_{\delta}, y_{\delta}\in X$ so that
\begin{equation}\label{eq1}
\varepsilon_0\leq d(x_{\delta}, y_{\delta})<\varepsilon_0+\delta
\end{equation}
and \begin{equation}\label{eq2}
d(Tx_{\delta}, Ty_{\delta})\geq\varepsilon_0.
\end{equation}
We take $\displaystyle \delta=\frac{1}{n}$, $n\geq 1$ natural number.
It follows that there are two sequences $\{x_n\}, \{y_n\}\subset X$ so that
\begin{equation}\label{eq3}
d(x_n,y_n)\geq\varepsilon_0,\quad n\geq 1~\textrm{and}~\lim_{n\rightarrow\infty} d(x_n,y_n)=\varepsilon_0
\end{equation}
and \begin{equation}\label{eq4}
d(Tx_n, Ty_n)\geq\varepsilon_0.
\end{equation}
Since $T$ is a $\mathcal{Z}-$contraction, there is $\zeta$ a simulation function so that
\begin{equation}\label{eq5}
0\leq\zeta(d(Tx_n,Ty_n),d(x_n,y_n)).
\end{equation}
From (\ref{eq4}), (\ref{eq5})and the condition $ii)$ in Definition \ref{defsimfct} it follows
$$0\leq\zeta(d(Tx_n,Ty_n),d(x_n,y_n))<d(x_n,y_n)-d(Tx_n,Ty_n)\leq d(x_n,y_n)-\varepsilon_0.$$
It follows \begin{equation}\label{eq6}
\lim_{n\rightarrow\infty}\zeta(d(Tx_n, Ty_n),d(x_n,y_n))=0
\end{equation}
But $\varepsilon_0\leq d(Tx_n, Ty_n)<d(x_n,y_n)-\zeta(d(Tx_n,Ty_n),d(x_n,y_n))$ witch implies $\displaystyle \lim_{n\rightarrow\infty} d(Tx_n,Ty_n)=\varepsilon_0.$

From condition $iii)$ in Definition \ref{defsimfct} we have $$\limsup_{n\rightarrow\infty}\zeta(d(Tx_n,Ty_n),d(x_n,y_n))<0$$ in contradiction with relation (\ref{eq6}).

The next result generalizes the main results from \cite{Dutta} and \cite{Rhoades}. See \cite{Aydi}.
\begin{theorem}
Let $(X,d)$ be a complete metric space and $T:X\rightarrow X$ be such that \begin{equation}
\psi(d(Tx,Ty))\leq\alpha(d(x,y))-\beta(d(x,y)),
\end{equation} for all $x,y\in X$, where $\psi, \alpha, \beta:[0,\infty)\rightarrow[0,\infty)$
are such that $\psi$ is continuous and non-decreasing, $\alpha$ is continuous, $\beta$ is lower semi-continuous,
\begin{equation}\psi(t)=0~\textrm{if and only if}~t=0,\alpha(0)=\beta(0)=0,
\end{equation}and\begin{equation}\psi(t)-\alpha(t)+\beta(t)>0~\textrm{for all}~t>0.
\end{equation}Then $T$ has a unique fixed point. \end{theorem}
We consider the following generalization of contractions.
\begin{df}
Let $(X,d)$ be a complete metric space and $T:X\rightarrow X$ be such that $$
\psi(d(Tx,Ty))\leq\alpha(d(x,y))-\beta(d(x,y)),
$$ for all $x,y\in X$, where $\psi, \alpha, \beta:[0,\infty)\rightarrow[0,\infty)$
are such that $\psi$ is non-decreasing, $\alpha$ is continuous, $\beta$ is lower semi-continuous,
and $$\psi(t)-\alpha(t)+\beta(t)>0~\textrm{for all}~t>0.
$$Then we say that $T$ is a weakly type contraction.
\end{df}
\begin{theorem}\label{thmMK}
If $T$ is a weakly type contraction, then $T$ is a Meir-Keeler contraction.
\end{theorem}
\bp We suppose that $T$ is not a Meir-Keeler contraction. Then there exists $\varepsilon_0>0$ and two sequences $\{x_n\},\{y_n\}\subset X$ such that $$\varepsilon_0\leq d(x_n,y_n)<\varepsilon_0+\frac{1}{n}$$ and $$d(Tx_n, Ty_n)\geq\varepsilon_0,\quad n\geq 1.$$
We have $\Psi(d(Tx_n,Ty_n))\geq\Psi(\varepsilon_0)$ and $\Psi(d(Tx_n,Ty_n))\leq\alpha(d(x_n,y_n))-\beta(d(x_n,y_n))$.
It follows $$\Psi(\varepsilon_0)\leq\alpha(d(x_n,y_n))-\beta(d(x_n,y_n))$$ or $$\beta(d(x_n,y_n))\leq\alpha(d(x_n,y_n))-\Psi(\varepsilon_0).$$

Using the continuity of $\alpha$, we have $$\liminf_{n\rightarrow\infty}\beta(d(x_n,y_n))\leq\alpha(\varepsilon_0)-\Psi(\varepsilon_0)$$
Since $\beta$ lower semi-continuous, it follows $\beta(\varepsilon_0)\leq\alpha(\varepsilon_0)-\Psi(\varepsilon_0),$ in contradiction with the hypothesis.
By Theorem \ref{thmufp} and Theorem \ref{thmMK}, we have the following Corollary:
\begin{corollary}Let $(X,d)$ be a complete metric space and $T$ be a weakly type contraction. Then $T$ has a unique fixed point $u$ and the Picard sequence $x_n=Tx_{n-1}$ converges to the fixed point of $T$.
\end{corollary}

\begin{flushright}
L. G\u avru\c ta, P.  G\u avru\c ta\\
\textit{
{ \normalsize Department of Mathematics, Politehnica University of Timi\c{s}oara, }\\
{ \normalsize Pia\c{t}a Victoriei no.2, 300006 Timi\c{s}oara, Rom\^{a}nia
}}\\

\hspace{5mm}{ \normalsize \textbf{E-mail}: gavruta\_laura@yahoo.com\\pgavruta@yahoo.com}
\end{flushright}

\begin{flushright}
F. Khojasteh\\
\textit{
{ \normalsize Department of Mathematics, Arak-Branch, }\\
{ \normalsize Islamic Azad University, Arak, Iran.
}}\\

\hspace{5mm}{ \normalsize \textbf{E-mail}: f-khojaste@iau-arak.ac.ir\\fr\_khojasteh@yahoo.com}
\end{flushright}

\end{document}